\documentclass[11pt]{article}
\usepackage{amsmath,amssymb,eucal}
\usepackage[dvips]{graphicx}

\newtheorem{thm}{\bf Theorem}[section]

\newtheorem{lm}[thm]{\bf Lemma}
\newtheorem{rem}[thm]{\bf Remark}
\newtheorem{prop}[thm]{\bf Proposition}

\numberwithin{equation}{section}

\begin{document}

\title{{\bf 
{\large Seifert surgery on knots via Reidemeister torsion and Casson-Walker-Lescop invariant II}}}
\author{
{\normalsize Teruhisa Kadokami, Noriko Maruyama and Tsuyoshi Sakai}}
\date{{\normalsize June 3, 2015}}
\footnotetext[0]{%
}
\maketitle

\begin{abstract}{
For a knot $K$ with $\Delta_K(t)\doteq t^2-3t+1$ in a homology $3$-sphere,
let $M$ be the result of $2/q$-surgery on $K$.
We show that an appropriate assumption on the Reidemeister torsion
of the universal abelian covering of $M$ implies $q=\pm 1$, if 
$M$ is a Seifert fibered space.

}\end{abstract}


\section{Introduction}\label{sec:intro}

The first auther \cite{Kd1} studied the Reidemister torsion of Seifert fibered homology lens spaces, and showed the following:
\begin{thm}\label{th:Seifert}
{\rm (\cite[Theorem 1.4]{Kd1})}
Let $K$ be a knot in a homology $3$-sphere $\Sigma$ such that
the Alexander polynomial of $K$ is $t^2-3t+1$.
The only surgeries on $K$ that may produce a Seifert fibered space with base $S^2$
and with $H_1\ne \{0\}, \mathbb{Z}$ have coefficients $2/q$ and $3/q$,
and produce Seifert fibered space with three singular fibers.
Moreover
(1) if the coefficient is $2/q$, then the set of multiplicities is
$\{2\alpha, 2\beta, 5\}$ where $\gcd(\alpha, \beta)=1$, and
(2) if the coefficient is $3/q$, then the set of multiplicities is
$\{3\alpha, 3\beta, 4\}$ where $\gcd(\alpha, \beta)=1$.
\end{thm}

\noindent
It is conjectured that Seifert surgeries on non-trivial knots are integral (except some cases). We \cite {KMS} have studied 
the $2/q$-Seifert surgery,
one of the remaining cases of the above theorem, by applying the Reidemister torsion and the Casson-Walker-Lescop invariant, 
and have given sufficient conditions to determine the integrality of $2/q$ (\cite[Theorems 2.1, 2.3]{KMS}).

In this paper, we give another condition for the integrality of  $2/q$ (Theorem 2.1). Like as in \cite{KMS}, 
the condition is also suggested by computations for the figure eight knot (\cite[Example 2.2]{KMS}).

We note two differences of this paper from \cite{KMS}; one is that the surgery coefficient appears in the condition instead of the Casson-Walker-Lescop invariant, and another is that we need more delicate estimation for the Dedekind sum to prove the result.

\medskip

\noindent
(1) Let $\Sigma$ be a homology $3$-sphere, 
and let $K$ be a knot in $\Sigma$.
Then $\Delta_K(t)$ denotes the Alexander polynomial of $K$, and
$\Sigma(K; p/r)$ denotes the result of $p/r$-surgery on $K$.

\bigskip
\noindent
(2) The first author \cite{Kd2} introduced the norm of polynomials and homology lens spaces: Let $\zeta_d$ be a primitive $d$-th root of unity.
For an element $\alpha$ of $\mathbb{Q}(\zeta_d)$,
$N_d(\alpha)$ denotes the norm of $\alpha$ associated to
the algebraic extension $\mathbb{Q}(\zeta_d)$ over $\mathbb{Q}$
.
Let $f(t)$ be a Laurent polynomial over $\mathbb{Z}$.
We define $|f(t)|_d$ by
$$|f(t)|_d=|N_d(f(\zeta_d))|
=\left| \prod_{i\in (\mathbb{Z}/d\mathbb{Z})^{\times}}
f(\zeta_d^i)\right|.$$
Let $X$ be a homology lens space with $H_1(X)\cong \mathbb{Z}/p\mathbb{Z}$.
Then there exists a knot $K$ in a homology $3$-sphere $\Sigma$ such that
$X=\Sigma(K; p/r)$ (\cite[Lemma 2.1]{BL}).
We define $|X|_d$ by
$$|X|_d=|\Delta_K(t)|_d,$$
where $d$ is a divisor of $p$.
Then $|X|_d$ is a topological invariant of $X$ (Refer to \cite{Kd2} for details).

\bigskip

\noindent
(3) Let $X$ be a closed oriented $3$-manifold.
Then $\lambda(X)$ denotes the Lescop invariant of $X$ (\cite{Le}). Note that $\lambda\left(S^3\right)=0$.

\section{Result}\label{sec:result}

Let $K$ be a knot in a homology $3$-sphere $\Sigma$.
Let $M$ be the result of $2/q$-surgery on $K$: $M=\Sigma(K; 2/q)$.
Let $\pi : X\to M$ be the universal abelian covering of $M$
(i.e.\ the covering associated to $\mathrm{Ker}(\pi_1(M)\to H_1(M))$).
Since $H_1(M)\cong \mathbb{Z}/2\mathbb{Z}$, $\pi$ is the $2$-fold unbranched covering.

In \cite{KMS}, we have defined $|K|_{(q,d)}$ by the following formula, if $|X|_d$ is defined:
$$|K|_{(q,d)}:=|X|_d.$$
Assume that the Alexander polynomial of $K$ is $t^2-3t+1$. Then, as noted in \cite{KMS}, $H_1(X)\cong \mathbb{Z}/5\mathbb{Z}$ and $|K|_{(q,5)}$ is defined.

\medskip
We then have the following.

\begin{thm}\label{th:main}
Let $K$ be a knot in a homology 3-sphere $\Sigma$. We assume the following.
 
\medskip
\noindent
{\bf (2.1)} $\lambda(\Sigma)=0$,

\medskip
\noindent
{\bf (2.2)} $\Delta_K(t)\doteq t^2-3t+1$,

\medskip
\noindent
{\bf (2.3)} $|q|\geq 3$,
 
\medskip
\noindent
{\bf (2.4)} $\sqrt{|K|_{(q,5)}}> 4q^2$.
\medskip

Then $M=\Sigma(K; 2/q)$ is not a Seifert fibered space.
\end{thm}

\begin{rem}
{\rm Let $K$ be the figure eight knot in $S^3$. Note that $\Delta_K(t)\doteq t^2-3t+1$. Then $|K|_{(q,5)}=(5q^2-1)^2$ by \cite[Example 2.2]{KMS}. 
Hence {\bf (2.4)} holds if $|q|\geq 3$.}
\end{rem} 
\begin{rem}{\rm 
Theorem 2.1 seems to suggest studying the asymptotic behavior of $|K|_{(q,d)}$ as a function of $q$.}
\end{rem}
 
\section{An inequality for the Dedekind sum}\label{sec:dedekind}
 
To prove Theorem 2.1, we need the following inequality for the Dedekind sum $s(\cdot, \cdot)$ ([RG]):
\begin{prop}\label{prop:mar}
{\rm (\cite[Lemma 3]{Ma})} For an even integer $p\geq 8$ and for an odd integer $q$ such that $3\leq q\leq p-3$ and $\gcd(p,q)=1$, we have
$$|s(q,p)|< f(2,p)$$
where $\displaystyle{f(2,p)=\frac{(p-1)(p-5)}{24p}}$.
\end{prop}
\medskip

By this proposition, we immediately have the following.

\begin{lm}\label{lm:ev}
For an even integer $p\geq 8$ and for an integer $q_{\ast}$ such that $q_{\ast} \not \equiv \pm 1 \;(\mathrm{mod}\ \! p)$ and $\gcd(p,q_{\ast})=1$, we have
$$|s(q_{\ast},p)|< \frac{p}{24}.$$
\end{lm}

\noindent
{\bf Proof.} By assumptions, there exists $q$ such that $q_{*} \equiv q \;(\mathrm{mod}\ \! p)$ and $3\leq q\leq p-3$. Hence by Proposition \ref{prop:mar}, we have
$$|s(q_{\ast},p)|=|s(q,p)|< \frac{(p-1)(p-5)}{24p}< \frac{p}{24}.$$
\hfill$\Box$

\begin{rem} {\rm The estimation given in Proposition 3.1 has a natural application (\cite{Ma}).}
\end{rem}

\section{Proof of Theorem \ref{th:main}}\label{sec:proof}

Suppose that $M=\Sigma (K;2/q)$ is a Seifert fibered space. Then, as shown in \cite{KMS}, we may assume that
\begin{center}
$(\ast)$ : $M$ has a framed link presentation as in Figure 1, 
\end{center}
where $1\le \alpha<\beta$ and $\gcd(\alpha, \beta)=1$.

\vspace{0.5cm}
\begin{figure}[htbp]
\begin{center}
\includegraphics[scale=0.6]{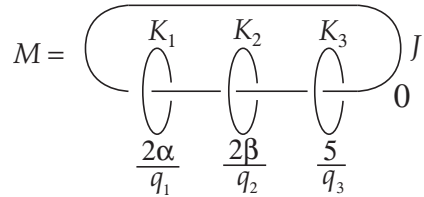}
\label{M1}
\caption{A framed link presentation of  $M=\Sigma (K;2/q)$}
\end{center}
\end{figure}

\medskip
Also as shown in \cite{KMS}, 
$\sqrt{|K|_{(q,5)}}=(\alpha \beta)^2$. Hence by {\bf (2.4)},
\begin{equation}\label{eq:ineq1}
(\alpha \beta)^2> 4q^2
\end{equation}

By {\bf (2.1)}, {\bf (2.2)} and \cite[1.5 T2]{Le}, we have $\lambda(M)=-q$. Hence $(\alpha \beta)^2> 4\{\lambda(M)\}^2$, and hence
\begin{equation}\label{eq:ineq2}
|\lambda(M)|< \frac{\alpha \beta}{2}
\end{equation}

\medskip
We now consider $e$ defined as follows:
$$e:=\frac{q_1}{2\alpha}+\frac{q_2}{2\beta}+\frac{q_3}{5}.$$ 
 
According to the sign of $e$, we treat two cases separetely: We first consider the case $e>0$.
Then the order of $H_1(M)$ is $20\alpha \beta e$. 
Since $H_1(M)\cong \mathbb{Z}/2\mathbb{Z}$, 
$20\alpha \beta e=2$, and
$e=1/(10\alpha \beta)$.
Hence by $(\ast)$ and \cite[Proposition 6.1.1]{Le}, we have
\begin{equation}\label{eq:lambda}
{\displaystyle
\lambda(M)=\left(-\frac{4}{5}\right)\alpha \beta
+\frac{5\beta}{24\alpha}
+\frac{5\alpha}{24\beta}
+\frac{1}{120\alpha \beta}-\frac 14-T
}
\end{equation}
where $T=s(q_1, 2\alpha)+s(q_2, 2\beta)+s(q_3, 5)$.

By (\ref{eq:ineq2}), we have
$$-\frac{\alpha \beta}{2}<\lambda(M).$$
Hence by  (\ref{eq:lambda}),
$$-\frac{\alpha \beta}{2}<
\left(-\frac{4}{5}\right)\alpha \beta+\frac{5\beta}{24\alpha}+\frac{5\alpha}{24\beta}
+\frac{1}{120\alpha \beta}-\frac 14+|T|.$$

Consequently
\begin{equation}\label{eq:ineq3}
{\displaystyle
\frac{3}{10}\alpha \beta<
-\frac{1}{4}+\frac{5}{24\alpha}\beta+\frac{5}{24}\left(\frac{\alpha}{\beta}\right)
+\frac{1}{120\alpha \beta}+ |T|
}
\end{equation}

As in  \cite{KMS}, we show that $\alpha \ge 2$ implies a contradiction: 
Suppose that $\alpha \ge 2$.
Since $\alpha<\beta$, we have $\beta \ge 3$ and $\alpha/\beta<1$.
Hence 
$$\frac{3}{5}\beta<-\frac 14+\frac{5}{24\cdot 2}\beta+\frac{5}{24}+\frac{1}{120\cdot 2\cdot 3}
+ |T|.$$
Since $|s(q_1,2\alpha)| \leq \frac {2\alpha} {12} < \frac {2\beta} {12}$, $|s(q_2,2\beta)| \leq \frac {2\beta} {12}$, 
and $|s(q_3,5)| \leq \frac 15$ as in  \cite{KMS}, we have  
 
$$|T|\le |s(q_1, 2\alpha)|+|s(q_2, 2\beta)|+|s(q_3, 5)|
\le \frac{\beta}{3}+\frac{1}{5}.$$
Hence
$$\frac{3}{5}\beta<-\frac 14+\frac{5}{48}\beta+\frac{5}{24}+\frac{1}{120\cdot 6}+\left(\frac \beta {3}+\frac 15\right).$$
Thus
$$\left( \frac 35 -\frac 5 {48} -\frac 13\right)\beta < -\frac 14+\frac{5}{24}+\frac{1}{120\cdot 6}+\frac 15.$$
Therefore
$$\frac{39}{240}\beta < \frac{1}{240}\left(38+\frac 13\right)<\frac{39}{240}.$$
This contradicts $\beta \geq 3$.

\medskip
We next show that $\alpha =1$ implies a contradiction:
Suppose that $\alpha=1$. By (\ref{eq:ineq1}), $\beta^2 > 4q^2$. Since $|q|\geq 3$, $\beta^2>4 \cdot 3^2=36$. Hence $\beta>6$. 
Since $\alpha=1$, $e=\displaystyle{\frac{1}{10\beta}}$. Hence
$$\frac{q_1}{2}+\frac{q_2}{2\beta}+ \frac{q_3}{5}=\frac{1}{10\beta}$$
and hence we have the following equation.
\begin{equation}\label{eq:eq4}
(5\beta)q_1+5q_2+(2\beta)q_3=1
\end{equation}
Since $q_1$ and $q_2$ are odd (see Figure 1), $\beta$ must be even. Since $\beta>6$, we have $\beta \geq 8$. We then have 
$$(\sharp): \;\;\; q_2 \not \equiv \pm 1 \;(\mathrm{mod}\ \! 2\beta).$$
In fact, since $q_1$ is odd, $(5\beta)q_1 \equiv \beta \;(\mathrm{mod}\ \! 2\beta)$. Hence by (\ref{eq:eq4}),
$$\beta+5q_2 \equiv 1 \;(\mathrm{mod}\ \! 2\beta).$$
Now suppose that $q_2 \equiv 1 \;(\mathrm{mod}\ \! 2\beta)$. Then $\beta+5 \equiv 1 \;(\mathrm{mod}\ \! 2\beta)$. 
This is impossible since $\beta \geq 8$. Next suppose that $q_2 \equiv -1 \;(\mathrm{mod}\ \! 2\beta)$. Then $\beta-5 \equiv 1 \;(\mathrm{mod}\ \! 2\beta)$. 
This is also impossible since $\beta \geq 8$. Thus $(\sharp)$ holds.

\medskip
Substituing $\alpha=1$ in 
(\ref{eq:ineq3}), 
$$\frac 3{10}\beta<-\frac 14+\frac{5}{24}\beta+\frac{5}{24\beta}
+\frac{1}{120\beta}+|T|$$
where $T=s(q_2, 2\beta)+s(q_3, 5)$
(since $s(q_1, 2)=0$).
By $(\sharp)$ and Lemma 3.2,
$$|s(q_2,2\beta)|<\frac{2\beta}{24}=\frac{\beta}{12}.$$
Hence
$$|T|\leq |s(q_2, 2\beta)|+|s(q_3, 5)|<\frac{\beta}{12}+\frac 15.$$
Since $\beta \geq 8$,
$$\frac 3{10}\beta<-\frac 14+\frac{5}{24}\beta+\frac{5}{24\cdot 8}
+\frac{1}{120\cdot 8}+\left(\frac{\beta}{12}+\frac 15\right).$$
Thus
$$\left(\frac 3{10}-\frac{5}{24}-\frac{1}{12}\right)\beta<-\frac 14+\frac{5}{24\cdot 8}
+\frac{1}{120\cdot 8}+\frac 15$$
and hence $\displaystyle{\frac 1{120}\beta<0}$. This is a contradiction, and ends the proof in the case $e>0$.

\medskip

We finally consider the case $e<0$.
Then ${\displaystyle e=-\frac{1}{10\alpha \beta}}$.
By $(\ast)$ and \cite[Proposition 6.1.1]{Le}, we have
$$\lambda(M)=-\left\{ \left(-\frac 45\right)\alpha \beta
+\frac{5\beta}{24\alpha}+\frac{5\alpha}{24\beta}
+\frac{1}{120\alpha \beta}-\frac 14+ T
\right\}.$$
Remaining part of the proof is similar to that in the case $e>0$.

\medskip
This completes the proof of Theorem 2.1.\hfill$\Box$


{\small
\par
Teruhisa Kadokami\par
Department of Mathematics, East China Normal University,\par
Dongchuan-lu 500, Shanghai, 200241, China\par
{\tt mshj@math.ecnu.edu.cn, kadokami2007@yahoo.co.jp}\par

\medskip

Noriko Maruyama\par
Musashino Art University,\par 
Ogawa 1-736, Kodaira, Tokyo 187-8505, Japan \par 
{\tt maruyama@musabi.ac.jp} \par

\medskip

Tsuyoshi Sakai\par
Department of Mathematics, Nihon University,\par
3-25-40, Sakurajosui, Setagaya-ku, Tokyo 156-8550, Japan \par
}

\end{document}